
\magnification=1200

\hsize=14.5cm \hoffset -0,8cm
\vsize=20cm
\baselineskip 15pt

\parindent=0pt   \parskip=0pt
\pageno=1


\hoffset=0mm    
\voffset=0cm    


\ifnum\mag=\magstep1
\hoffset=0cm   
\voffset=0cm   
\fi


\pretolerance=500 \tolerance=1000  \brokenpenalty=5000

\catcode`\@=11

\font\eightrm=cmr8         \font\eighti=cmmi8
\font\eightsy=cmsy8        \font\eightbf=cmbx8
\font\eighttt=cmtt8        \font\eightit=cmti8
\font\eightsl=cmsl8        \font\sixrm=cmr6
\font\sixi=cmmi6           \font\sixsy=cmsy6
\font\sixbf=cmbx6


\font\tengoth=eufm10       \font\tenbb=msym8 at 10pt
\font\eightgoth=eufm8      \font\eightbb=msym8
\font\sevengoth=eufm7      \font\sevenbb=msym7
\font\sixgoth=eufm6        \font\fivegoth=eufm5


\font\tengo=eufm10 \font\sevengo=eufm7 \font\fivego=eufm5
\newfam\gofam \textfont\gofam=\tengo \scriptfont\gofam=\sevengo
\scriptscriptfont\gofam=\fivego 

\skewchar\eighti='177 \skewchar\sixi='177
\skewchar\eightsy='60 \skewchar\sixsy='60


\newfam\gothfam           \newfam\bbfam

\def\tenpoint{%
  \textfont0=\tenrm \scriptfont0=\sevenrm \scriptscriptfont0=\fiverm
  \def\rm{\fam\z@\tenrm}%
  \textfont1=\teni  \scriptfont1=\seveni  \scriptscriptfont1=\fivei
  \def\oldstyle{\fam\@ne\teni}\let\old=\oldstyle
  \textfont2=\tensy \scriptfont2=\sevensy \scriptscriptfont2=\fivesy
  \textfont\gothfam=\tengoth \scriptfont\gothfam=\sevengoth
  \scriptscriptfont\gothfam=\fivegoth
  \def\goth{\fam\gothfam\tengoth}%
  \textfont\bbfam=\tenbb \scriptfont\bbfam=\sevenbb
  \scriptscriptfont\bbfam=\sevenbb
  \def\bb{\fam\bbfam\tenbb}%
  \textfont\itfam=\tenit
  \def\it{\fam\itfam\tenit}%
  \textfont\slfam=\tensl
  \def\sl{\fam\slfam\tensl}%
  \textfont\bffam=\tenbf \scriptfont\bffam=\sevenbf
  \scriptscriptfont\bffam=\fivebf
  \def\bf{\fam\bffam\tenbf}%
  \textfont\ttfam=\tentt
  \def\tt{\fam\ttfam\tentt}%
  \abovedisplayskip=12pt plus 3pt minus 9pt
  \belowdisplayskip=\abovedisplayskip
  \abovedisplayshortskip=0pt plus 3pt
  \belowdisplayshortskip=4pt plus 3pt
  \smallskipamount=3pt plus 1pt minus 1pt
  \medskipamount=6pt plus 2pt minus 2pt
  \bigskipamount=12pt plus 4pt minus 4pt
  \normalbaselineskip=12pt
  \setbox\strutbox=\hbox{\vrule height8.5pt depth3.5pt width0pt}%
  \let\bigf@nt=\tenrm       \let\smallf@nt=\sevenrm
  \normalbaselines\rm}

\def\eightpoint{%
  \textfont0=\eightrm \scriptfont0=\sixrm \scriptscriptfont0=\fiverm
  \def\rm{\fam\z@\eightrm}%
  \textfont1=\eighti  \scriptfont1=\sixi  \scriptscriptfont1=\fivei
  \def\oldstyle{\fam\@ne\eighti}\let\old=\oldstyle
  \textfont2=\eightsy \scriptfont2=\sixsy \scriptscriptfont2=\fivesy
  \textfont\gothfam=\eightgoth \scriptfont\gothfam=\sixgoth
  \scriptscriptfont\gothfam=\fivegoth
  \def\goth{\fam\gothfam\eightgoth}%
  \textfont\bbfam=\eightbb \scriptfont\bbfam=\sevenbb
  \scriptscriptfont\bbfam=\sevenbb
  \def\bb{\fam\bbfam}%
  \textfont\itfam=\eightit
  \def\it{\fam\itfam\eightit}%
  \textfont\slfam=\eightsl
  \def\sl{\fam\slfam\eightsl}%
  \textfont\bffam=\eightbf \scriptfont\bffam=\sixbf
  \scriptscriptfont\bffam=\fivebf
  \def\bf{\fam\bffam\eightbf}%
  \textfont\ttfam=\eighttt
  \def\tt{\fam\ttfam\eighttt}%
  \abovedisplayskip=9pt plus 3pt minus 9pt
  \belowdisplayskip=\abovedisplayskip
  \abovedisplayshortskip=0pt plus 3pt
  \belowdisplayshortskip=3pt plus 3pt
  \smallskipamount=2pt plus 1pt minus 1pt
  \medskipamount=4pt plus 2pt minus 1pt
  \bigskipamount=9pt plus 3pt minus 3pt
  \normalbaselineskip=9pt
  \setbox\strutbox=\hbox{\vrule height7pt depth2pt width0pt}%
  \let\bigf@nt=\eightrm     \let\smallf@nt=\sixrm
  \normalbaselines\rm}

\tenpoint


\def\pc#1{\bigf@nt#1\smallf@nt}         \def\pd#1 {{\pc#1} }


\catcode`\;=\active
\def;{\relax\ifhmode\ifdim\lastskip>\z@\unskip\fi
\kern\fontdimen2  -1.2 \fontdimen3 \string;}

\catcode`\:=\active
\def:{\relax\ifhmode\ifdim\lastskip>\z@\unskip\fi\penalty\@M\ \fi\string:}

\catcode`\!=\active
\def!{\relax\ifhmode\ifdim\lastskip>\z@
\unskip\fi\kern\fontdimen2  -1.1 \fontdimen3 \string!}

\catcode`\?=\active
\def?{\relax\ifhmode\ifdim\lastskip>\z@
\unskip\fi\kern\fontdimen2  -1.1 \fontdimen3 \string?}

\def\^#1{\if#1i{\accent"5E\i}\else{\accent"5E #1}\fi}
\def\"#1{\if#1i{\accent"7F\i}\else{\accent"7F #1}\fi}

\frenchspacing


\newtoks\auteurcourant      \auteurcourant={\hfil}
\newtoks\titrecourant       \titrecourant={\hfil}

\newtoks\hautpagetitre      \hautpagetitre={\hfil}
\newtoks\baspagetitre       \baspagetitre={\hfil}

\newtoks\hautpagegauche
\newtoks\hautpagedroite

\newtoks\baspagegauche      \baspagegauche={\centerline{\folio}}
\newtoks\baspagedroite      \baspagedroite={\centerline{\folio}}

\newif\ifpagetitre          \pagetitretrue


\headline={\ifpagetitre\the\hautpagetitre
\else\ifodd\pageno\the\hautpagedroite\else\the\hautpagegauche\fi\fi}

\footline={\ifpagetitre\the\baspagetitre\else
\ifodd\pageno\the\baspagedroite\else\the\baspagegauche\fi\fi
\global\pagetitrefalse}


\def\raggedbottom{\topskip 10pt plus 36pt\r@ggedbottomtrue}



\def\pointir{\unskip . --- \ignorespaces}


\def\Bigbreak{\vskip-\lastskip\bigbreak}
\def\Medbreak{\vskip-\lastskip\medbreak}


\def\ctexte#1\endctexte{%
  \hbox{$\vcenter{\halign{\hfill##\hfill\crcr#1\crcr}}$}}


\long\def\ctitre#1\endctitre{%
    \ifdim\lastskip<24pt\vskip-\lastskip\bigbreak\bigbreak\fi
  		\vbox{\parindent=0pt\leftskip=0pt plus 1fill
          \rightskip=\leftskip
          \parfillskip=0pt\bf#1\par}
    \bigskip\nobreak}

\long\def\section#1\endsection{%
\vskip 0pt plus 3\normalbaselineskip
\penalty-250
\vskip 0pt plus -3\normalbaselineskip
\Bigbreak
\message{[section \string: #1]}{\bf#1\unskip}\pointir}

\long\def\sectiona#1\endsection{%
\vskip 0pt plus 3\normalbaselineskip
\penalty-250
\vskip 0pt plus -3\normalbaselineskip
\Bigbreak
\message{[sectiona \string: #1]}%
{\bf#1}\medskip\nobreak}

\long\def\subsection#1\endsubsection{%
\Medbreak
{\it#1\unskip}\pointir}

\long\def\subsectiona#1\endsubsection{%
\Medbreak
{\it#1}\par\nobreak}

\def\rem#1\endrem{%
\Medbreak
{\it#1\unskip} : }

\def\remp#1\endrem{%
\Medbreak
{\pc #1\unskip}\pointir}

\def\rema#1\endrem{%
\Medbreak
{\it #1}\par\nobreak}

\def\newparwithcolon#1\endnewparwithcolon{
\Medbreak
{#1\unskip} : }

\def\newparwithpointir#1\endnewparwithpointir{
\Medbreak
{#1\unskip}\pointir}

\def\newpara#1\endnewpar{
\Medbreak
{#1\unskip}\smallskip\nobreak}


\long\def\th#1 #2\enonce#3\endth{%
   \Medbreak
   {\pc#1} {#2\unskip}\pointir{\it #3}\medskip}

\long\def\tha#1 #2\enonce#3\endth{%
   \Medbreak
   {\pc#1} {#2\unskip}\par\nobreak{\it #3}\medskip}


\long\def\Th#1 #2 #3\enonce#4\endth{%
   \Medbreak
   #1 {\pc#2} {#3\unskip}\pointir{\it #4}\medskip}

\long\def\Tha#1 #2 #3\enonce#4\endth{%
   \Medbreak
   #1 {\pc#2} #3\par\nobreak{\it #4}\medskip}


\def\decale#1{\smallbreak\hskip 28pt\llap{#1}\kern 5pt}
\def\decaledecale#1{\smallbreak\hskip 34pt\llap{#1}\kern 5pt}
\def\puce{\smallbreak\hskip 6pt{$\scriptstyle\bullet$}\kern 5pt}



\def\displaylinesno#1{\displ@y\halign{
\hbox to\displaywidth{$\@lign\hfil\displaystyle##\hfil$}&
\llap{$##$}\crcr#1\crcr}}


\def\ldisplaylinesno#1{\displ@y\halign{
\hbox to\displaywidth{$\@lign\hfil\displaystyle##\hfil$}&
\kern-\displaywidth\rlap{$##$}\tabskip\displaywidth\crcr#1\crcr}}


\def\eqalign#1{\null\,\vcenter{\openup\jot\m@th\ialign{
\strut\hfil$\displaystyle{##}$&$\displaystyle{{}##}$\hfil
&&\quad\strut\hfil$\displaystyle{##}$&$\displaystyle{{}##}$\hfil
\crcr#1\crcr}}\,}


\def\system#1{\left\{\null\,\vcenter{\openup1\jot\m@th
\ialign{\strut$##$&\hfil$##$&$##$\hfil&&
        \enskip$##$\enskip&\hfil$##$&$##$\hfil\crcr#1\crcr}}\right.}


\let\@ldmessage=\message

\def\message#1{{\def\pc{\string\pc\space}%
                \def\'{\string'}\def\`{\string`}%
                \def\^{\string^}\def\"{\string"}%
                \@ldmessage{#1}}}



\def\up#1{\raise 1ex\hbox{\smallf@nt#1}}


\def\qed{\raise -2pt\hbox{\vrule\vbox to 10pt{\hrule width 4pt
                 \vfill\hrule}\vrule}}

\def\virg{\raise .4ex\hbox{,}}   


\def\build#1_#2^#3{\mathrel{
\mathop{\kern 0pt#1}\limits_{#2}^{#3}}}


\def\boxit#1#2{%
\setbox1=\hbox{\kern#1{#2}\kern#1}%
\dimen1=\ht1 \advance\dimen1 by #1 \dimen2=\dp1 \advance\dimen2 by #1
\setbox1=\hbox{\vrule height\dimen1 depth\dimen2\box1\vrule}%
\setbox1=\vbox{\hrule\box1\hrule}%
\advance\dimen1 by .4pt \ht1=\dimen1
\advance\dimen2 by .4pt \dp1=\dimen2  \box1\relax}

\def\ra{\rightarrow}

\def\blanc{\hskip 1em plus 0,25em minus 0,25em}
\def\fhd#1#2{\nospacedmath\smash{\mathop{\hbox to 12mm{\rightarrowfill}}
\limits^{\scriptstyle#1}_{\scriptstyle#2}}}
\def\fhg#1#2{\nospacedmath\smash{\mathop{\hbox to 12mm{\leftarrowfill}}
\limits^{\scriptstyle#1}_{\scriptstyle#2}}}


\def\P{{\bb P}}

\def\mv=msxm10\def\dra{{\mv\char 16}}
\mathcode`A="7041 \mathcode`B="7042 \mathcode`C="7043 \mathcode`D="7044
\mathcode`E="7045 \mathcode`F="7046 \mathcode`G="7047 \mathcode`H="7048
\mathcode`I="7049 \mathcode`J="704A \mathcode`K="704B \mathcode`L="704C
\mathcode`M="704D \mathcode`N="704E \mathcode`O="704F \mathcode`P="7050
\mathcode`Q="7051 \mathcode`R="7052 \mathcode`S="7053 \mathcode`T="7054
\mathcode`U="7055 \mathcode`V="7056 \mathcode`W="7057 \mathcode`X="7058
\mathcode`Y="7059 \mathcode`Z="705A

\def\spacedmath#1{\def\packedmath##1${\bgroup\mathsurround=0pt ##1\egroup$}%
\mathsurround#1 \everymath={\packedmath}\everydisplay={\mathsurround=0pt }}

\def\nospacedmath{\mathsurround=0pt \everymath={}\everydisplay={} }



\spacedmath{2pt}

\parindent=0mm


\catcode`\@=12

\showboxbreadth=-1  \showboxdepth=-1



\def\Z{{\bf Z}}\def\Q{{\bf Q}}\def\C{{\bf C}}
\def\O{{\cal O}}
\def\dem{{}$\sqcap\!\!\!\!\sqcup$}\def\bl{\bigl}\def\Bl{\Bigl}\def\br{\bigr}
\def\Br{\Bigr}\def\fin{\dem\medskip}\def\deb{\dem}
\vsize 20,5cm

\def\Ext{{\rm Ext}}
\def\End{{\rm End}}
\def\Hom{{\rm Hom}}
\def\Pic{{\rm Pic}}
\def\E{{\bf E}}
\def\e{{\cal E}\{l\}}
\def\g{{\cal E}[l]}
\def\M{{{\cal F}({\cal E},l)}}

\def\gr{gr^l}
\def\D{D[l]}
\def\d{{\cal D}^l}
\def\Pic{{\rm Pic}}
\def\Ker{{\rm Ker}}
\def\ext{{\cal E}xt}

\centerline{A PROPOS DE L'EXISTENCE DE FIBR\'ES STABLES SUR LES SURFACES}
\bigskip
\centerline{par Andr\'e HIRSCHOWITZ
\footnote{*}{\sevenrm Dans le cadre du programme de recherche {\it Fibr\'es
vectoriels} de
Europroj.}
et Yves LASZLO
\footnote{**}{\sevenrm Partiellement financ\'e par
le contrat Science {\it Geometry of Algebraic Varieties}, contrat n\up{o}
SCI-0398-C (A).}.}
\bigskip
{\pc R\'ESUM\'E:} On prouve ici qu'il existe (en caract\'eristique nulle) sur
toute
surface  polaris\'ee (lisse, projective et connexe) des fibr\'es stables de
rang $r\geq 2$, de
premi\`ere classe de Chern arbitraire et de $c_2$ assez grand.
\bigskip
{\pc ABSTRACT:} We prove the existence (in characteristic 0) on every
polarized (smooth, projective and connected) surface of stable bundles of rank
$r\geq 2$, arbitrary
first Chern class and large enough $c_2$.
\bigskip
\section I. Introduction
\endsection
Ce travail est la suite logique de l'\'etude des faisceaux prioritaires sur
$\P^2$ men\'ee dans
[H\&L]. On s'int\'eresse maintenant \`a des faisceaux sur une surface
quelconque. Le but de cette
note est de donner une preuve tr\`es simple (et alg\'ebrique) de l'existence de
fibr\'es stables de
rang et premi\`ere classe de Chern quelconques et de seconde classe de Chern
assez grande sur toute
surface lisse (connexe et projective). En rang 2, on a ([Gi\& Li ]) des
r\'esultats bien plus
pr\'ecis, \`a savoir que les espaces de modules correspondants sont
irr\'eductibles (pour $c_2$ assez
grand). Maruyama a quant \`a lui prouv\'e l'existence de fibr\'es stables de
$c_2$ arbitrairement
grand en rang quelconques et premi\`ere classe de Chern [M].

L'existence de fibr\'es stables de premi\`ere classe de Chern nulle, de rang
quelconque et de $c_2$
assez grand a \'et\'e obtenue par divers auteurs ([A], [Gi], [T], [F]...). Par
ailleurs, tout
r\'ecemment Gieseker et Li ont
annonc\'e
des r\'esultats  tr\`es pr\'ecis en rang et premi\`ere
classe de Chern quelconques et
 seconde
classe de Chern assez grande, \`a savoir un th\'eor\`eme
``\`a la Donaldson-Friedman-Zhu"
 d'une part, et d'autre part l'irr\'eductibilit\'e de ces espaces de
modules, comme dans le cas du rang 2 .

Notre
 m\'ethode
consiste \`a
utiliser les techniques infinit\'esimales d\'evelopp\'ees dans
[D\&P] pour montrer qu'une d\'eformation g\'en\'erique d'une diminution
\'el\'ementaire assez
singuli\`ere d'un fibr\'e fix\'e au d\'epart est stable. De plus, ce faisceau
stable est localement
libre (si $c_2$ est assez grand).

Contrairement \`a [H\&L], ces r\'esultats n'ont \'et\'e obtenus qu'en
caract\'eristique nulle. Outre les arguments de type lissit\'e g\'en\'erique,
le point essentiel est
qu'un th\'eor\`eme (effectif) de restriction des faisceaux stables \`a des
courbes (cf. II) du type
[F] ne semble pas \^etre disponible, \`a supposer qu'il soit vrai, en
caract\'eristique positive. Si
on disposait d'un tel \'enonc\'e en caract\'eristique positive, on obtiendrait
alors l'existence de
fibr\'es {\it semi-stables} de rang et premi\`ere classe de Chern quelconques
et de seconde classe de
Chern assez grande sur toute surface lisse (connexe et projective). Signalons
que le th\'eor\`eme
principal n'est pas effectif, sauf lorsque le groupe de N\'eron-Severi de la
surface est $\Z$.

Enfin, notre r\'esultat a \'et\'e  \'egalement
obtenu, par une m\'ethode tr\`es diff\'erente,
par Zhenbo
Qin et Wei Ping Li dans [Q\&Li].

\section II. Notations et lemmes pr\'eliminaires.
\endsection
Soit $X$ une surface lisse, projective et connexe sur $\C$. Dans toute la
suite, la premi\`ere
classe de Chern d'un faisceau
\footnote{(${}^1$)}{\sevenrm Par faisceau, on entend faisceau coh\'erent.}
est vue comme un \'el\'ement du groupe de N\'eron-Severi
$NS(X)$ et la seconde comme un rationnel.  On fixe \'egalement un fibr\'e
tr\`es
ample $\O(1)$  sur $X$ de premi\`ere classe de Chern $h$: il d\'efinit une
notion de degr\'e et de
stabilit\'e (au sens de Mumford). Soit $E$ (resp. $E'$) un faisceau sur $X$ de
rang $r$ (resp. $r'$). Comme dans
[D\&P], le discriminant $\Delta(E)$ est d\'efini par la formule
$$\Delta(E)={1\over r}(c_2(E)-{c_1(E)^2\over 2}+{c_1(E)^2\over 2r}).$$

\rem Remarque
\endrem Le discriminant est invariant par twist par un fibr\'e inversible. En
particulier, cette
d\'efinition permet d'\'ecrire agr\'eablement la formule de Riemann-Roch (voir
{\it infra}).

\th D\'EFINITION 1
\enonce On d\'efinit la pente d'un faisceau  $E$ par $$\mu(E)={c_1(E)\over
r}\in NS(X)\otimes\Q.$$
\endth

En appliquant la formule de Riemann-Roch dans le groupe de Grothendieck de $X$
\`a la classe
$[E'].[E]^\vee$, on montre qu'il existe un certain polyn\^ome $Q$ de
$Num(X)[T',T]$ ne d\'ependant que
de $X$ v\'erifiant $$\chi(E,E')=\sum_{i=0}^2(-1)^i\dim
\Ext^i(E,E')=rr'\Bl(Q\bigl(\mu(E'),\mu(E)\bigr)-\Delta(E)-\Delta(E')\Br).$$

\th PROPOSITION 2
\enonce Il existe un polyn\^ome $P\in\Z[T',T]$ ne d\'ependant que de $(X,h)$
tel que
pour tout faisceau semi-stable $E$ de rang $r$ et de pente $\mu$ on ait:
$$H^0(X,E)\leq P(r,(\mu.h)).$$
\endth
\deb On peut supposer $(\mu.h)\geq 0$. Le faisceau $E(-n),\ n>(\mu.h)$ n'a pas
de section.
D'apr\`es un r\'esultat de Flenner ([F], th\'eor\`eme 1.2), la restriction de
$E$ \`a un diviseur
g\'en\'erique de $\mid\!\! \O(n)\!\!\mid,\ n\geq {r^2\over 2}\deg(X)$ est
semi-stable. Soit par
exemple $n=1+r(\mu.h)+r^2\deg(X)$. Soit $H$ un  diviseur de $\mid\!\!
\O(n)\!\!\mid$ tel que $H$
soit lisse et connexe et $E\mid_H$ semi-stable. La suite de cohomologie
associ\'ee \`a la suite
exacte  $$0\ra E(-n)\ra E\ra E\!\mid_H\ra 0$$
donne la majoration
$$h^0(E)\leq h^0(E\!\mid_H).\leqno (1)$$
On sait ([L]], prop. 3.10)
$$2.h^0(E\!\mid_H)\leq (r^2-1)(g-1)+r\mu(E\!\mid_H)+1+\dim
\End(E\!\mid_H),\leqno (2)$$
$g$ \'etant le genre de $H$. Comme $E\!\mid_H$ est semi-stable, on a
$$\dim \End(E\!\mid_H)\leq r^2,\leqno (3)$$
(utiliser la filtration de Jordan-H\"older de $E\!\mid_H$), de sorte que la
formule  du genre jointe \`a
(1), (2) et (3) donne le r\'esultat.
\fin

On obtient alors le
\th MAJORATION 3
\enonce Pour tous faisceaux semi-stables $E,E'$ de rang $r,r'$ et de pente
$\mu,\mu'$ on a:
$$\dim \Hom(E,E')\leq P(rr',(\mu'-\mu).h).$$
\endth

\section III. La construction
\endsection
Soit $L$ un faisceau inversible sur $X$. On fixe un entier $r\geq 2$
et un fibr\'e ${\cal E}$ de rang $r\geq 2$ et de premi\`ere classe de Chern
$c_1$ arbitraire (par
exemple le fibr\'e $(r-1).\O_X\oplus L$ avec $c_1(L)=c_1$).

\th D\'EFINITIONS 1
\enonce (i) On appelle diminution \'el\'ementaire d'un fibr\'e $E$, tout
sous-faisceau $E'$ noyau d'un
morphisme surjectif $\phi:\ E\ra \O_Z$ o\`u $Z$ est un sous-sch\'ema de $X$ de
dimension $0$.

(ii) Pour tout entier $l\geq 0$, on d\'efinit le faisceau ${\cal E}\{l\}$,
diminution
\'el\'ementaire g\'en\'erique de $\cal E$, comme le noyau du morphisme
g\'en\'erique
${\cal E}\ra \O_Z$, o\`u $Z$ est le sous-sch\'ema lisse g\'en\'erique de
longueur $l$ de $X$.
 \endth

On va montrer que ${\cal E}\{l\}$ est un point lisse du champ des faisceaux
coh\'erents. Notons tout
d'abord ce lemme, dont la d\'emonstration est laiss\'ee au lecteur

\th LEMME 2
\enonce Soit $E$ un fibr\'e et $M$ un fibr\'e en droites. Soit $f\in
\Hom(E,E\otimes M)$ un
endomorphisme qui n'est pas une homot\'ethie. Alors, il existe un point $z\in
X$ et une surjection
$\phi:\ E\ra \O_z$ telle que $f$ n'induise pas un endomorphisme de $\Ker(\phi
)$.
\endth

\rem Notation
\endrem Pour tout entier $i$, on note
$$\Ext^i({\cal F},{\cal F})\up{o}\buildrel\hbox{\sevenrm d\'ef}\over
=\Ker\bl\{\Ext^i({\cal F},{\cal
F})\ra H^i(\O_X)\br\}.$$

En utilisant le lemme et la dualit\'e de Serre, on obtient le

\th COROLLAIRE 3
\enonce Si $l\geq \dim \Ext^2({\cal E},{\cal E})\up{o}$ le groupe
$\Ext^2({\cal E}\{l\},{\cal E}\{l\})\up{o}$
est nul. \endth

\rem Notation
\endrem On note $l_0$ le plus petit entier $l$, tel qu'on ait
$$\Ext^2({\cal E}\{m\},{\cal E}\{m\})\up{o}=0\hbox{ pour }m\geq l.$$
\medskip
\centerline{\bf On suppose d\'esormais $\bf l\geq l_0$.}
\medskip
On d\'eduit que le champ des faisceaux de d\'eterminant fix\'e est lisse en
${\cal
E}\{l\}$. Comme  le sch\'ema de Picard est lisse (en caract\'eristique nulle),
il existe une unique
composante $\M$ du
 champ des faisceaux (sans condition de d\'eterminant) passant par ${\cal
E}\{l\}$.

\th D\'EFINITION 4
\enonce Pour $l\geq l_0$, on d\'efinit le faisceau ${\cal E}[l]$
comme le faisceau repr\'esent\'e par le point g\'en\'erique de $\M$.
\endth

\rem Remarque
\endrem On commet l\`a un abus de langage: par d\'efinition du point
g\'en\'erique, le faisceau
${\cal E}[l]$ est une classe d'\'equivalence de faisceaux d\'efinis sur des
corps qui dominent $\M$.

\section IV Lissit\'e du fibr\'e g\'en\'erique ${\cal E}[l]$
\endsection
A tout fibr\'e $\cal E$ et tout entier $l\geq 0$, on a donc associ\'e une
unique composante
irr\'eductible $\M$ du champ des faisceaux de rang $r$, de premi\`ere classe de
Chern $c_1(L)$ et de
seconde classe de Chern $c_2({\cal E})+l$. Cette composante a la dimension
attendue, \`a savoir
$$\dim Ext^1(\e,\e)\up{o}+\dim\Pic(X)-1.$$
\rem Remarques
\endrem

1. La dimension $\dim Ext^1(\e,\e)\up{o}$ se calcule par Riemann-Roch.

2. Le fait que la dimension attendue soit $\dim
Ext^1(\e,\e)\up{o}+\dim\Pic(X)-1$ et non
$\dim Ext^1(\e,\e)\up{o}+\dim\Pic(X)$ vient du fait qu'on regarde le champ des
fibr\'es et non
les espaces de modules associ\'es (c'est le m\^eme ph\'enom\`ene qui se produit
d\'ej\`a pour les
fibr\'es en droites, puisque le champ des fibr\'es en droites a pour dimension
$\dim Pic(X)-1$).

\th LEMME 1
\enonce
Soit ${\cal E}\{l\}$ la diminution \'el\'ementaire g\'en\'erique de $\cal E$.
Si $l$
est assez grand, la fl\`eche naturelle
$$ev:\ \Ext^1({\cal E}\{l\},{\cal E}\{l\})\otimes\O\rightarrow
{\cal E}xt^1({\cal E}\{l\},{\cal E}\{l\})$$
est surjective.
\endth
\medskip
\dem La fl\`eche $H^0(ev)$ est un edge-homomorphisme de la suite
spectrale $$E_2^{p,q}=H^p\bigl(\ext^q({\cal E}\{l\},{\cal
E}\{l\})\bigr)\Longrightarrow
\Ext^{p+q}({\cal E}\{l\},{\cal E}\{l\}).$$
Soit $d$ la diff\'erentielle
$$d_2^{0,1}: H^0\bl(\ext^1({\cal E}\{l\},{\cal E}\{l\})\br)\ra
H^2\bl(\ext^0({\cal E}\{l\},{\cal E}\{l\})\br).$$

On a
$$E_\infty^{0,1}=\Ker(d)$$
de sorte qu'on a une surjection
$$\Ext^1({\cal E}\{l\},{\cal E}\{l\})\ra \Ker(d).\leqno (1)$$
La codimension de $\Ker(d)$ dans $H^0\bl(\ext^1({\cal E}\{l\},{\cal
E}\{l\})\br)$ est donc major\'ee
par $$\dim H^2\bl({\cal E}xt^0({\cal E}\{l\},{\cal E}\{l\})\br)=\dim
H^2\bl({\cal
E}xt^0({\cal E},{\cal E})\br),$$ dimension ind\'ependante de $l$. La fl\`eche
$ev$ est donc surjective
en dehors d'un ensemble de cardinal au plus $\dim H^2\bl({\cal E}xt^0({\cal
E},{\cal E})\br)$.
Ainsi, d\`es que $l$ est $>\dim H^2\bl({\cal E}xt^0({\cal E},{\cal E})\br)$, il
existe un point singulier de
${\cal E}\{l\}$ (d\'efini sur un corps assez gros) o\`u $ev$ est surjective.
\medskip
Montrons que $ev$ est en fait partout surjective (tous les points jouent le
m\^eme r\^ole).

Au dessus de $S^l_*\P {\cal E}\times X$ existe une famille universelle $\E$ de
diminutions
\'el\'ementaires ($S^l_*\P {\cal E}$ est le produit sym\'etrique moins toutes
les diagonales). Soit
$U\subset S^l_*\P {\cal E}\times X$ le support (muni de sa structure r\'eduite)
du conoyau $\E\ra
\E^{\vee\vee}$.
Le sch\'ema $U$ est irr\'eductible et la premi\`ere projection $pr_1:U\ra
S^l_*\P {\cal E}$ est finie.
D'apr\`es ce qui pr\'ec\`ede, il existe un point $u\in U$ tel que $ev(u)$ est
surjective. Le ferm\'e
des points de $U$ tels que $ev(u)$ n'est pas surjectif est $\not = U$ et comme
$pr_1$ est fini, son
image par $pr_1$ est $\not = S^l_*\P {\cal E}$. C'est ce qu'on voulait. On
aurait pu aussi dire
que le rev\^etement $pr_1$ est galoisien (de groupe $S_l$) et faire agir le
groupe de Galois.\fin

\th PROPOSITION 2
\enonce Avec les notations pr\'ec\'edentes, si $l$ est assez grand, le fibr\'e
g\'en\'erique ${\cal E}[l]$ est localement libre.
\endth
\deb Examinons la situation locale: soit $\O$ l'anneau local d'un point de $Z$
et ${\cal I}$ son
id\'eal. Localement, $\e$ est une somme directe $(r-1)\O\oplus{\cal I}$. La
d\'eformation
semi-universelle locale $\cal U$ de $(r-1)\O\oplus{\cal I}$ est lisse (le
groupe $\ext^2$
correspondant est nul car $(r-1)\O\oplus{\cal I}$ est de dimension projective
$1$). Comme $X$ est
lisse et que le support de $\e$ est r\'eduit, le complexe de Koszul
$$0\ra\O\ra\O\oplus\O\ra{\cal
I}\ra 0$$ est exact de sorte que $r\O$ se sp\'ecialise sur $(r-1)\O\oplus{\cal
I}$. On d\'eduit que
$(r-1)\O\oplus{\cal I}$ se g\'en\'eralise en un module (localement) libre.

Passons du local au global: d'apr\`es le lemme pr\'ec\'edent, l'application
tangente correspondante
$$\Ext^1(\e,\e)\otimes\O\ra {\cal E}xt^1(\e,\e)$$
est surjective d\`es que $l$ est assez grand. Le morphisme du global vers le
local est donc une
submersion au point d\'efini par $\e$. On d\'eduit alors de l'\'etude locale
que le point
g\'en\'erique de $\M$ est localement libre.\fin

\section V Le fibr\'e g\'en\'erique est stable.
\endsection En consid\'erant une filtration de Jordan-H\"older de chaque
gardu\'e de la
filtration de Harder-Narasimhan de ${\cal E}[l]$, on obtient une filtration
dite de
Jordan-Harder-Narasimhan de ${\cal E}[l]$ qui d\'efinit un drapeau
$$D[l]=\Bigl(0=F_0[l]\subset
F_1[l]\subset F_2[l]\ldots\subset F_n[l]={\cal E}[l]\Bigr),$$ tel que les
gradu\'es
$\gr_i=F_i[l]/F_{i-1}[l]$ soient semi-stables. Soient
$$H_i(n)=r_i^l\bl(P(\mu_i^l+nh)-\Delta_i^l\br),\ i=1,\ldots n$$  les
polyn\^omes de Hilbert des
gradu\'es $\gr_i$. Par d\'efinition, la suite $\bl((\mu_i^l.h)\br)_{i=1,\ldots
,n}$ est
d\'ecroissante (au sens large). \rem Remarque et notation
\endrem L'entier $n$ n'est \'egal \`a $1$ que si $\e$ est stable. On le note
$n(l)$.
\medskip
Comme dans [D\&P], il
existe un champ alg\'ebrique au dessus de $\M$ param\'etrant les familles de
drapeaux param\'etr\'ees
par des bases $S$ au dessus de $\M$ dont les gradu\'es sont $S$-plats de
polyn\^omes de Hilbert
$H_i$. Soit $f$ le morphisme (repr\'esentable) d'oubli de la filtration
$$f:\d\ra \M.$$
Le drapeau $\D$ d\'efinit un point $\bar\eta$ de $\d$ au dessus de $\eta$,
point
g\'en\'erique de $\M$
\footnote{(${}^1$)}{\sevenrm En fait, il faudrait faire une extension des
scalaires, les
filtrations de Jordan-H\"older d'un faisceau semi-stable $F$ n'\'etant en
g\'en\'eral d\'efinis que
sur une extension finie du corps de d\'efinition de $F$. Comme on travaille en
caract\'eristique
nulle, l'extension de corps est \'etale de sorte que l'\'etude infinit\'esimale
n'est aucunement
affect\'ee par cette extension.}.

\medskip
La structure infinit\'esimale de $\d$ a \'et\'e \'etudi\'ee en d\'etail dans
[D\&P]. On d\'eduit de
la proposition 1.5 de [D\&P] et de la surjectivit\'e de l'application de
Kodaira-Spencer au point
g\'en\'erique (adapter la proposition 1.1 de [H\&L]), l'existence une suite
exacte
$$T_{\bar\eta}\d{\nospacedmath\build\hbox to
8mm{\rightarrowfill}_{}^{T_{\bar\eta}f}}\Ext^1\bigl(\g,\g\bigr)\hbox to
8mm{\rightarrowfill}
\Ext^1_+\bigl(\g,\g\bigr)\hbox to 8mm{\rightarrowfill}
\Ext^2_-\bigl(\g,\g\bigr),$$ les espaces
vectoriels $\Ext^1_+\bigl(\g,\g\bigr)$ et $\Ext^2_-\bigl(\g,\g\bigr)$ \'etant
des
groupes d'extensions filtr\'ees (relativement \`a la filtration
pr\'ec\'edente).

D'apr\`es la proposition III.10.6 de [Ha], le morphisme $f$ \'etant dominant
par construction,
$T_{\bar\eta}f$ est surjective de sorte qu'on doit avoir une injection
$$\Ext^1_+\bigl(\g,\g\bigr)\hookrightarrow
\Ext^2_-\bigl(\g,\g\bigr).\leqno (*)$$  On va
voir que c'est impossible d\`es que $\limsup\limits_{l\ra\infty}n(l)\geq 2$.

\rem Notations.
\endrem
Soit $\cal F$ un faisceau. Soit $\mu_{max}$ (resp. $\mu_{min}$) la pente du
premier (resp. dernier)
gradu\'e de la fitration de Jordan-Harder-Narasimhan de $\cal F$. On notera
$\delta({\cal F})$ le
rationnel $(\mu_{max}-\mu_{min}).h$.
\medskip
On veut contr\^oler la taille des caract\'eristiques  d'Euler associ\'ees aux
$\gr_i(\eta)$ en
fonction des seuls discriminants. Montrons d'abord le

\th LEMME 1
\enonce Soit $F$ un faisceau et $r\in\Q$ un rationnel donn\'es. La famille des
pentes $\mu$ des
sous-faisceaux $G$ de $F$ tels que $\mu.h\geq r$ est finie.
\endth le des
satur\'es
$$\hat G={\rm \Ker}\Bl\{F\ra (F/G)_{/Torsion}\Br\}$$
est {\it limit\'ee} de sorte que l'ensemble des $\mu(G)$ est {\it fini}. Il
existe un diviseur
effectif $D_G$ tel que
$$c_1(\hat G)=c_1(G)+[D_G],$$
la classe $[D_G]$ de $D_G$ dans $NS(X)$ \'etant simplement
$c_1\bl(Torsion(F/G)\br)$.
Comme $(\mu(\hat G).h)$ est born\'e, l'hypoth\`ese assure que le degr\'e de
$D_G$ est born\'e.
L'existence des vari\'et\'es de Chow assure alors que la famille des $D_G$ est
limit\'ee. Le lemme
suit.\fin

Montrons alors la proposition suivante:

\th PROPOSITION 2
\enonce Il existe un ensemble fini qui contient les pentes $\mu_i^l$ des
gradu\'es $\gr_i$, ce pour tout $l\geq 0$.
\endth

\deb Par d\'efinition de la filtration de Jordan-Harder-Narasimhan, on a
$$(\mu_i^l.h)\geq\bl(\mu(\g).h\br)=\bl(\mu({\cal E}).h\br).$$

Tout drapeau $\D$ se sp\'ecialise en
un drapeau $D_s^l$ de $\e$, et les pentes des sous-faisceaux correspondants
sont constantes par
sp\'ecialisation (platitude des gradu\'es). Comme $\e\subset {\cal E}$ et que
ces deux faisceaux
co\"\i ncident en codimension 1, la proposition est maintenant cons\'equence du
lemme
3.\fin

On doit maintenant contr\^oler les discriminants $\Delta_i^l$.

\th LEMME 3
\enonce Soit $E$ un faisceau semi-stable de rang $r$ et de pente $\mu$. Il
existe une constante $C_0$
ne d\'ependant que de $(X,h)$ et de $(r,\mu)$ telle que $$\Delta(E)\geq C_0.$$
\endth
\deb La formule de Riemann-Roch donne l'in\'egalit\'e
$$r\Delta(E)\geq r.P(\mu)-\dim H^0(E)-\dim H^2(E).$$
La MAJORATION II.3 donne alors la minoration cherch\'ee.\fin

On est maintenant en mesure de contr\^oler les groupes
$\Ext^1_+\bigl(\g,\g\bigr)$ et $\Ext^2_-\bigl(\g,\g\bigr)$ intervenant dans
(*).

\th PROPOSITION 4
\enonce (i) Il existe des constante $C,\ l_1$ ne
d\'ependant que de $(X,h)$ et $\cal E$ telles que,
$$\dim \Ext^2_-\bigl(\g,\g\bigr)\leq C\ pour\ l\geq l_1.$$

(ii) Si $\limsup\limits_{l\ra\infty}n(l)\geq 2$, on a l'\'egalit\'e
$$\limsup_{l\ra\infty}\dim \Ext^1_+\bigl(\g,\g\bigr)=+\infty.$$
\endth

\deb  On \'etudie ces groupes d'extension gr\^ace \`a
l'existence de suites spectrales (voir [D\&P]) $Sp_+$ (resp. $Sp_-$)
convergeant vers
$\Ext_+^i(E,E)$ (resp. $\Ext_-^i(E,E))$. Le terme $E_1$ de $Sp_+$ (resp.
$Sp_-$)
est donn\'e par

$$E_1^{p,q}=\left\{\matrix{
\prod\limits_i \Ext^{p+q}(\gr_i(E),\gr_{i-p}(E))&si\blanc p<0,\cr
0&sinon.\cr
}\right. \leqno (Sp_+)$$
(resp.)
$$E_1^{p,q}=\left\{\matrix{
0&si\blanc p<0,\cr
\prod\limits_i \Ext^{p+q}(\gr_i(E),\gr_{i-p}(E))&sinon.\cr
}\right. \leqno (Sp_-)$$
correspondants.

Preuve de (i): on doit majorer les dimensions des groupes
$$\Ext^2(\gr_i,\gr_j),\ i\leq j.$$ Par dualit\'e de Serre, on s'int\'eresse aux
groupes
$$\Hom(\gr_j,\gr_i\otimes \omega_X),\ i\leq j.$$
Or,
$$k-\delta(\g)\leq (\mu_j^l.h)-(\mu_i^l.h)+k<k,$$
$k$ \'etant le degr\'e de $\omega_X$. Comme $\delta$ croit par
sp\'ecialisation, on a finalement
$$k-\delta({\cal E})\leq(\mu_j^l.h)-(\mu_i^l.h)+k<k,$$
ce qui prouve (i) gr\^ace \`a la MAJORATION II.3 et la suite spectrale $Sp_-$.

Preuve (ii): d\'efinissons la caract\'eristique d'Euler filtr\'ee $\chi_+({\cal
F},{\cal F})$ par
$$\chi_+({\cal F},{\cal F})\buildrel\hbox{\sevenrm d\'ef}\over
=\sum_{i=0}^2\dim Ext^i_+({\cal F},{\cal F}).$$
L'invariance de la caract\'eristique d'Euler dans la suite spectrale $Sp_+$ et
Riemann-Roch donnent
$$\chi_+(\e,\e)=\sum_{i>j}\chi\bl(\gr_i,\gr_j\br)=\sum_{i>j}r_ir_j
\Bl(Q\bl(\mu_i^l,\mu_j^l\br)-\Delta_i^l-\Delta_j^l\Br).\leqno (1)$$ De plus,
$$\chi\bl(\g\br)=\system{ \sum_{i=1}^n\chi\bl(\gr_i\br)=\sum_{i=1}^n
r_i\Bl(Q\bl(\mu_i^l\br)-\Delta_i^l\Br)\cr \chi({\cal E})-l.\hfill\cr}\leqno
(2)$$

\medskip
Dela proposition V.2, du lemme V.3 et de (2), on d\'eduit
$$\lim_{l\ra +\infty}\sup_{i}\Delta_i^l=+\infty.\leqno (3)$$
Dela proposition V.2, du lemme V.3 et des formules (1) et (3), on d\'eduit
alors
$$\lim_{l\ra +\infty}\chi_+(\e,\e)=-\infty.$$\fin

\th COROLLAIRE 5
\enonce On a $\lim\limits_{l\ra\infty}n(l)=1$ et $\e$ est stable si $l$ est
assez grand.
\endth
\medskip

Utilisant le corollaire pr\'ec\'edent, on d\'eduit alors la

\th PROPOSITION 6
\enonce Soit $\cal E$ un fibr\'e vectoriel de rang $\geq 2$. Si $l$ est assez
grand, le faisceau $\g$
est localement libre et stable. De plus, $\M$ est g\'en\'eriquement lisse de
dimension attendue.
\endth

On a donc prouv\'e en particulier le

\th T\'EOR\`EME 7
\enonce Pour tout $r\geq 2$ et toute classe $c_1\in NS(X)$, il existe un entier
$n$ tel que si
$l\geq n$, il existe un fibr\'e stable de rang $r$, de premi\`ere classe de
Chern $c_1$ et de seconde
classe de Chern $l$ et qui d\'efinit un point lisse du champ des faisceaux sur
$X$.
\endth

\rem Remarque
\endrem
Il n'y a aucune difficult\'e \`a rendre effectif le th\'eor\`eme dans le cas
o\`u le groupe de
N\'eron-Severi de $X$ est $\Z$. Toutefois, m\^eme dans ce cas, la borne obtenue
semble beaucoup trop
grande pour \^etre v\'eritablement int\'eressante.

\bigskip
\centerline{\bf BIBLIOGRAPHIE}
\bigskip

[A] I.V. Artamkin, {\it Deforming torsion-free sheaves on an algebraic
surface}, Math. USSR Izvestiya
{\bf 36} (1991), 449-485.
\medskip
[D\&P] J. M. Drezet \& J. Le Potier, {\it Fibr\'es stables et fibr\'es
exceptionnels sur $\P$}, Ann. scient. \'Ec. Norm. Sup., 4\up{o} s\'erie, {\bf
18} (1985), 193-244.
\medskip
[F] H. Flenner, {\it Restrictions of semi-stable bundles on projective
varieties}, Comment.
Math. Helvetici {\bf 59} (1984) p. 635-650.
\medskip
[Gi] D. Gieseker {\it A construction of stable bundles on an algebraic
surface}, J. Diff. Geom.
{\bf 27}, (1988), 137-154. \medskip
[Gi\& Li ] D. Gieseker \& Jun Li {\it Moduli of Vector Bundles over Surfaces
I}, preprint 1993.
\medskip
[G] A. Grothendieck, {\it Technique de descente et
th\'eor\`emes d'existence en g\'eom\'etrie alg\'ebrique. IV. Les sch\'emas de
Hilbert}, S\'em.
Bourbaki {\bf 221} (1960).
\medskip
[Ha] R. Hartshorne, {\it Algebraic Geometry}, GTM 52, Spinger Verlag,
Berlin-Heidelberg-New York,
(1977).
\medskip
[H\&L] A. Hirschowitz \& Y. Laszlo {\it A propos de l'espace des modules de
fibr\'es de rang
2 sur $\P$}, Math. Annal. {\bf 297} (1993), p. 85-102.
\medskip
[Q\&Li] Z. Qin \& Wei Ping Li, {Stable vector bundles over algebraic surfaces},
preprint 1993.
\medskip
[L] G. Laumon, {\it Fibr\'es vectoriels sp\'eciaux}, Bull. Soc. math. France,
{\bf 119} (1991),
p. 97-119.
\medskip
[Q\&Li] Z. Qin \& W. P. Li {\it Stable vector bundles on algebraic surface},
preprint 1993.
\medskip
[M] M Maruyama, {\it Moduli of stable sheaves, II}, J. of Math.
Kyoto Univ., 18 (1978), 557-614.
\medskip
[T] C. Taubes, {\it Self dual connections on 4-manifolds with infinite
intersection matrix}, J.
Diff. Geom {\bf 19} (1984), 517-560.

\footnote{}{\sevenrm Andr\'e HIRSCHOWITZ Universit\'e de Nice U.R.A 168 Parc
Valros 06108 Nice Cedex
02 France

Yves LASZLO Universit\'e Paris-Sud U.R.A 752 Math\'ematiques B\^atiment 425
91405 Orsay Cedex
France et Universit\'e Bordeaux 1 U.R.A. 226 Math\'ematiques 351 cours de la
Lib\'eration 33405
Bordeaux Cedex France }\end